\documentclass[10pt,reqno]{article}

\usepackage{amssymb}
\usepackage{epsfig}
\usepackage{amsmath}
\usepackage{amsthm}
\usepackage{color}
\definecolor{r}{rgb}{0.9,0.3,0.1}
\definecolor{b}{rgb}{0.1,0.3,0.9}

\newtheorem{theorem}{Theorem}[section]
\newtheorem{lemma}[theorem]{Lemma}

\theoremstyle{remark}
\newtheorem{remark}[theorem]{Remark}

\theoremstyle{definition}
\newtheorem{assumption}[theorem]{Assumption}

\newtheorem{definition}[theorem]{Definition}

\newcommand\cbrk{\text{$]$\kern-.15em$]$}}
\newcommand\opar{\text{\,\raise.2ex\hbox{${\scriptstyle
|}$}\kern-.34em$($}}
\newcommand\cpar{\text{$)$\kern-.34em\raise.2ex\hbox{${\scriptstyle |}$}}\,}

\newcommand{\ga}{\gamma}

\newcommand{\om}{\omega}

\newcommand{\de}{\delta}

\newcommand\bL{\mathbb{L}}
\newcommand\bR{\mathbb{R}}
\newcommand\bH{\mathbb{H}}
\newcommand\bZ{\mathbb{Z}}

\newcommand\bE{\mathbb{E}}

\newcommand\cF{\mathcal{F}}
\newcommand\cH{\mathcal{H}}

\newcommand\cP{\mathcal{P}}

\newcommand\cO{\mathcal{O}}

\newcommand\frH{\mathfrak{H}}

\newcommand{\mysection}[1]{\section{#1}
\setcounter{equation}{0}}

\begin{document}
\setlength{\baselineskip}{16pt}

\title
{A $W^1_2$-theory of  Stochastic  Partial Differential Systems of Divergence type on $C^1$ domains}

\author{Kyeong-Hun Kim\footnote{Department of Mathematics, Korea
University, 1 Anam-dong, Sungbuk-gu, Seoul, South Korea 136-701,
\,\, kyeonghun@korea.ac.kr. The research of this author is supported
by the Korean Research Foundation Grant funded by the Korean
Government 20090087117}\quad \hbox{\rm and} \quad Kijung
Lee\footnote{Department of Mathematics, Ajou University, Suwon,
South Korea 443-749, \, kijung@ajou.ac.kr.} }


\date{}


\maketitle

\begin{abstract}
In this paper we  study the  stochastic  partial differential
systems of divergence type with $C^1$ space domains in $\bR^d$.
Existence and uniqueness results are obtained in terms of Sobolev
spaces with weights so that we allow the derivatives of the solution
to blow up near the boundary. The coefficients of the systems are only measurable and
are allowed to  blow up near the boundary.

\vspace*{.125in}


\vspace*{.125in}

\end{abstract}



\mysection{Introduction}

In this article we are dealing with $W^1_2$-theory of the stochastic partial differential systems (SPDSs)  of
$d_1$  equations of divergent type:
\begin{eqnarray}
du^k&=&(D_i(a^{ij}_{kr}u^r_{x^j}+\bar{b}^i_{kr}u^r+\bar{f}^{ik})+b^i_{kr}u^r_{x^i}+c_{kr}u^r+ f^k)dt\nonumber\\
&&+(\sigma^i_{kr,m}u^r_{x^i}+\nu_{kr,m}u^r+g^k_m)dw^m_t,  \quad \quad t>0
\label{eqn main system}\\
u^k(0)&=&u^k_0\nonumber
\end{eqnarray}
with $x\in \bR^d,\;\bR^d_+$ or   $\cO$, a bounded $C^1$ domain. Here,
$\{w^m_{t}:m=1,2,\ldots\}$ is a countable set of independent
one-dimensional Brownian motions defined on a probability space
$(\Omega,\cF,P)$. Indices $i$ and $j$ run from $1$ to $d$ while
$k,\;j=1,2,\cdots,d_1$ and $m=1,2,\cdots$. To make expressions
simple, we are using the summation convention on $i,j,r,m$. The
coefficients $a^{ij}_{kr}, \bar{b}^i_{kr},
b^i_{kr},c_{kr},\sigma^i_{kr,m}$ and $\nu_{kr,m}$ are measurable
functions depending on $\omega\in \Omega,t,x$. Detailed formulation
of (\ref{eqn main system}) follows in the subsequent sections.

Demand for a general theory of stochastic partial differential
systems(SPDSs) arises when we model the interactions among unknowns
in a natural phenomenon with random behavior. For example, the motion
of a random string can be modeled by means of SPDSs(see \cite{MT}
and \cite{Fu}).

We note that, if $d_1=1$, then the system (\ref{eqn main system}) becomes a
single stochastic partial differential equation (SPDE) of
divergence type. In this case $L_2$-theory on $\bR^d$ was developed long ago and an account of it can
be found, for instance, in \cite{R}. Also,
$L_p$-theory($p\geq 2$) of such single equations with $C^1$ space
domains can be found in \cite{Kim04}, \cite{Kim} and \cite{Yoo99}
in which weighted Sobolev spaces are used to allow derivatives of
the solutions to blow up near the boundary. For comparison with
$L_p$-theory of SPDEs of non-divergence type, we refer to
\cite{Kim03}, \cite{KK}, \cite{KL2}, \cite{Kr99} and references
therein.

The main goal of this article is to extend the results \cite{R},
\cite{Kim04}, \cite{Kim}, \cite{Yoo99} for single equations to the
case of systems under no smoothness assumptions on the coefficients.
We prove the uniqueness and existence results of system (\ref{eqn
main system}) in weighted Sobolev spaces so that we allow the
derivatives of the solutions to blow up near the boundary. The
coefficients of the system are only measurable and are allowed to blow up near the boundary (See
(\ref{05.04.01})).

We declare that $W^1_p$-theory, a desirable further result beyond
$W^1_2$-theory, is not successful yet even under the assumption that
the coefficients $a^{ij}_{kr}$ and $\sigma^{i}_{kr}$ are constants.
This is due to the difficulty caused by considering SPDSs instead
of SPDEs.  For $L_p$-theory, $p>2$, one must overcome
tremendous mathematical difficulties rising in the general settings;
one of the main difficulties in the case $p>2$ is that the arguments
we are using in the proof of  Lemma \ref{a priori 1}  below are not
working since in this case we get some extra terms which we simply can
not control.

The organization of the article is as follows. Section \ref{Cauchy}
handles the Cauchy problem. In section \ref{section half} and section \ref{section domain} we develop
our theory of the system defined on $\bR^d_+$ and bounded domain $\cO$, respectively.

 As usual, $\bR^{d}$
stands for the Euclidean space of points $x=(x^{1},...,x^{d})$,
$B_{r}(x)=\{y\in\bR^{d}:|x-y|<r\}$, $B_{r}=B_{r}(0)$ and
$\bR^{d}_{+}=\{x\in\bR^{d}:x^{1}>0\}$.
For $i=1,...,d$, multi-indices $\alpha=(\alpha_{1},...,\alpha_{d})$,
$\alpha_{i}\in\{0,1,2,...\}$, and functions $u(x)$ we set
$$
u_{x^{i}}=\frac{\partial u}{\partial x^{i}}=D_{i}u,\quad
D^{\alpha}u=D_{1}^{\alpha_{1}}\cdot...\cdot D^{\alpha_{d}}_{d}u,
\quad|\alpha|=\alpha_{1}+...+\alpha_{d}.
$$
If we write  $c=c(\cdots)$, this means that the constant $c$ depends
only on what are in parenthesis.

\mysection{The systems on $\bR^d$}\label{Cauchy}

In this section we develop  some solvability results of linear
systems defined on space domain $\bR^d$. These results will be used
later for systems defined on $\bR^d_+$ or a bounded $C^1$ domain
$\mathcal{O}$.

Let $(\Omega,\mathcal{F},P)$ be a complete probability space and
$\{\mathcal{F}_t\}$ be a filtration such that $\mathcal{F}_0$
 contains all $P $-null sets of $\Omega $; the probability space
 $(\Omega,\mathcal{F},P)$ is rich so that we define independent one-dimensional $\{\mathcal{F}_t\}$-adapted Wiener
processes $\{w^{m}_{t}\}_{m=1}^{\infty}$ on it. We let $\mathcal{P}$
denote the predictable $\sigma$-algebra on $\Omega\times(0,\infty)$.

The space $C^{\infty}_0=C^{\infty}_0(\mathbb{R}^d;\mathbb{R}^{d_1})$
denotes the set of all $\mathbb{R}^{d_1}$-valued infinitely
differentiable functions with compact support in $\mathbb{R}^d$. By
$\mathcal{D}$ we mean the space of $\mathbb{R}^{d_1}$-valued
distributions on $C^{\infty}_0$; precisely, for $u\in \mathcal{D}$
and $\phi\in C^{\infty}_0$ we define $(u,\phi)\in \mathbb{R}^{d_1}$
with components $(u,\phi)^k=(u^k,\phi^k)$, $k=1,\ldots,d_1$. Each
$u^k$ is a usual $\mathbb{R}$-valued distribution defined on
$C^{\infty}(\mathbb{R}^d;\mathbb{R})$.
We let $L_p=L_p(\mathbb{R}^d;\mathbb{R}^{d_1})$ be the space of all
$\mathbb{R}^{d_1}$-valued functions $u=(u^1,\ldots,u^{d_1})$
satisfying
\[
\|u\|^p_{L_p}:=\sum^{d_1}_{k=1}\|u^k\|^p_{L_p}<\infty.
\]
For $p \in[2,\infty)$ and $\gamma\in(-\infty,\infty)$ we define the
space of Bessel potential
$H^{\gamma}_p=H^{\gamma}_p(\mathbb{R}^d;\mathbb{R}^{d_1})$ as the
space of all distributions $u$ such that $(1-\Delta)^{n/2}u\in L_p$,
where 
$$
((1-\Delta)^{\gamma/2}u)^k:=(1-\Delta)^{\ga/2}u^k
:=\cF^{-1}[(1+|\xi|^2)^{\gamma/2}\cF(u^k)(\xi)].
$$
Here, $\cF$ is the Fourier transform. Define
\[
\|u\|_{H^{\gamma}_p}:=\|(1-\Delta)^{\ga/2}u\|_{L_p}.
\]
Then, $H^{\gamma}_p$ is a Banach space with the given norm and
$C^{\infty}_0$ is dense in $H^{\gamma}_p$. Note that $H^{\ga}_p$ are
usual Sobolev spaces for $\ga=0,1,2,\ldots$. It is well known that
the first order differentiation operators,
$\partial_i:H^{\gamma}_{p}(\mathbb{R}^d;\mathbb{R})\to
H^{\gamma-1}_p(\mathbb{R}^d;\mathbb{R})$ given by $u\to u_{x^i}$
$(i=1,2,\ldots,d)$, are bounded. On the other hand, for $u\in
H^{\gamma}_{p}(\mathbb{R}^d;\mathbb{R})$, if $\text{supp}\, (u)
\subset (a,b)\times \mathbb{R}^{d-1}$ with $-\infty<a<b<\infty$,
we have
\begin{equation}
                                        \label{eqn 5.1.1}
\|u\|_{H^{\gamma}_{p}(\mathbb{R}^d;\mathbb{R})}\leq
c(d,\gamma,a,b)\|u_{x}\|_{H^{\gamma-1}_{p}(\mathbb{R}^d;\mathbb{R})}
\end{equation}
(see, for instance, Remark 1.13 in \cite{kr99}). Let $\ell_2$
be the set of all real-valued sequences $e=(e_1,e_2,\ldots)$ with
the inner product $(e,f)_{\ell_2}=\sum_{m=1}^{\infty}e_mf_m$ and the
norm $|e|_{\ell_2}:=(e,e)^{1/2}_{\ell_2}$. For
$g=(g^1,g^2,\cdots,g^{d_1})$, where  $g^k$ are  $\ell_2$-valued functions,
 we define
$$
\|g\|^p_{H^{\gamma}_p(\ell_2)}:=\sum_{k=1}^{d_1}\||(1-\Delta)^{\ga/2}g^k|_{\ell_2}\|^p_{L_p}.
$$

Using the spaces mentioned above, for a fixed time $T$, we define
the stochastic Banach spaces
$$
\bH^{\ga}_p(T):=L_p(\Omega\times(0,T], \mathcal{P},H^{\ga}_p), \quad
\bH^{\ga}_p(T,\ell_2):=L_p(\Omega\times(0,T],
\mathcal{P},H^{\ga}_p(\ell_2)),
$$
$$
\mathbb{L}_p(T):=\bH^0_p(T),\quad
\mathbb{L}_p(T,\ell_2)=\bH^0_p(T,\ell_2)
$$
with norms given by
\[
\|u\|^p_{\mathbb{H}^{\ga}_p(T)}=\bE\int^{T}_0\|u(t)\|^p_{H^{\gamma}_p}dt,\quad
\|g\|^p_{\mathbb{H}^{\ga}_p(T,\ell_2)}=\bE\int^{T}_0\|u(t)\|^p_{H^{\gamma}_p(\ell_2)}dt.
\]
Lastly, we set
$U^{\gamma}_p=L_p(\Omega,\mathcal{F}_0,H^{\gamma-2/p}_p)$.

\begin{definition}\label{md}
For a $\mathcal{D}$-valued function $u\in\mathbb{H}^{\ga+2}_p(T)$,
we write $u\in\mathcal{H}^{\ga+2}_p(T)$ if $u(0,\cdot)\in
U^{\gamma+2}_p$ and there exist $f\in \mathbb{H}^{\ga}_p(T),\;g\in
\mathbb{H}^{\ga+1}_p(T,\ell_2)$ such that
$$
du=f\,dt+g^m dw^m_t, \quad t\leq T
$$
in the sense of distributions, that is, for  any $\phi\in
C^{\infty}_0$ and $k=1,2,\cdots,d_1$, the equality
\begin{equation}\label{e}
(u^k(t,\cdot),\phi)= (u^k(0,\cdot),\phi)+ \int^t_0( f^k(s,\cdot),\phi)ds+
\sum_{m=1}^{\infty}\int^t_0(g^k_{m}(s,\cdot),\phi)dw^{m}_s
\end{equation}
holds (a.s.) for all $t\leq T$.  The norm in
$\cH^{\gamma+2}_{p}(T)$ is defined by
\[
\|u\|_{\mathcal{H}^{\ga+2}_p(T)}= \|u\|_{\mathbb{H}^{\ga+2}_p(T)}+\|
f\|_{\mathbb{H}^{\ga}_p(T)}+\|g\|_{\mathbb{H}^{\ga+1}_p(T,\ell_2)}+
\|u(0)\|_{U^{\ga+2}_p}.
\]
\end{definition}

\begin{remark}
Note that since the  coefficients in system (\ref{eqn main system})
are only measurable, the space  $\cH^{\gamma+2}_p(T)$ is not
appropriate for    system (\ref{eqn main system}) unless
$\gamma=-1$.

\end{remark}

We set $A^{ij}=(a^{ij}_{kr})$, $\Sigma^i=(\sigma^i_{kr})$ and
$\mathcal{A}^{ij}=(\alpha^{ij}_{kr})$, where
$$
\alpha^{ij}_{kr}=\frac{1}{2}\sum_{l=1}^d
(\sigma^i_{lk},\sigma^j_{lr})_{\ell_2},\quad
\sigma^i_{kr}=(\sigma^i_{kr,1},\sigma^i_{kr,2},\cdots).
$$
Also, we set $\bar{B}^i=(\bar{b}^i_{kr}),
B^i=(b^i_{kr}),C=(c_{kr}),\mathcal{N}=(\nu_{kr})$, where
$\nu_{kr}:=(\nu_{kr,1},\nu_{kr,2},\ldots)$.

 For any $d_1\times d_1$ matrix $M=(m_{kr})$ we let
$$
|M|:=\sqrt{\sum_{k,r}(m_{kr})^2}\;;\quad
|M|:=\sqrt{\sum_{k,r}|m_{kr}|_{\ell_2}^2},
$$
where the latter is the case that the elements are in $\ell_2$.

 Throughout
the article we assume the following.

\begin{assumption}
                     \label{main assumptions}
(i) The coefficients $a^{ij}_{kr}, \bar{b}^i_{kr}, b^i_{kr},
c_{kr},\sigma^i_{kr,m}$ and $\nu_{kr,m}$ are $\mathcal{P}\times
\mathcal{B}(\mathbb{R}^d)$-measurable, where
$\mathcal{B}(\mathbb{R}^d)$ denotes Borel $\sigma$-field in
$\mathbb{R}^d$.

(ii) There exist finite constants $\delta,\;K^j(j=1,\ldots,d),\;
L>0$ so that
\begin{equation}\label{assumption 1}
\de|\xi|^2\le\xi^*_{i}\left(A^{ij}-\mathcal{A}^{ij}\right)\xi_{j}
\end{equation}
holds for any $\omega\in\Omega,\ t> 0$, where
 $\xi$ is any
(real) $d_1\times d$ matrix, $\xi_i$ is the $i$th column of $\xi$;
again the summations on $i,j$ are understood. Moreover, we assume
that for any $\om$, $t>0$, $x\in\mathbb{R}^d$, $i,j=1,\ldots,d$,
\begin{equation}\label{assumption 2}
\left|A^{1j}(\om,t,x)\right|\le K^j,  \quad
\left|A^{ij}(\om,t,x)\right|\le L\;(i\ne 1), \quad
|\mathcal{A}^{ij}(\om,t,x)|\le L.
\end{equation}

\end{assumption}


Our main theorem in this section is the following.
\begin{theorem}
                      \label{thm 1}
Assume that there is a constant $N_0\in (1,\infty)$ such that for
any $\om$, $t>0$, $x\in\mathbb{R}^d$, $i=1,\ldots,d$,
\begin{equation}
                   \label{extra}
|\bar{B}^i|,\; |B^i|,\; |C|,\;
|\mathcal{N}|<N_0.
\end{equation}
Then  for any $\bar{f}^i\in \bL_2(T)$ $(i=1,\ldots,d)$, $ f\in
\bH^{-1}_2(T)$, $g\in \bL_2(T,\ell_2)$, and $u_0\in U^{1}_2$, system
(\ref{eqn main system}) has a unique solution $u\in
\mathcal{H}^{1}_2(T)$, and for this solution we have
\begin{eqnarray}
\|u_x\|_{\bL_{2}(T)}&\leq&
c(\|u\|_{\bL_{2}(T)}+\sum_i\|\bar{f}^i\|_{\bL_2(T)}
+\|f\|_{\bH^{-1}_2(T)}+\|g\|_{\bL_2(T,\ell_2)}+\|u_0\|_{U^{1}_2}),\label{e 6.5.2}\\
 \|u\|_{\bH^{1}_{2}(T)}&\leq& ce^{cT}(\sum_i\|\bar{f}^i\|_{\bL_2(T)}+\|f\|_{\bH^{-1}_2(T)}
 +\|g\|_{\bL_2(T,\ell_2)}+\|u_0\|_{U^{1}_2}),\label{e 6.5.3}
\end{eqnarray}
where $c=c(d,d_1,\delta,K,L,N_0)$.
\end{theorem}

\begin{proof}

1. We note that $f^k$ can be expressed as $ f^k=F^{0k}+\mathrm{div}
(F^{1k},F^{2k},\ldots,F^{dk})$, where $F^{0k}\in
\mathbb{H}^{1}_2(T)$, $F^{ik}\in \mathbb{H}^{\gamma}_2(T)$ with the
estimate
$\|F^{0k}\|_{\mathbb{H}^{1}_2(T)}+\sum_{i=1}^{d}\|F^{ik}\|_{\mathbb{L}_2(T)}\le
c(d,d_1)\| f^k\|_{\mathbb{H}^{-1}_2(T)}$; this follows from the
observation $ f^k=(1-\Delta)(1-\Delta)^{-1} f^k=(1-\Delta)^{-1}
f^k+\mathrm{div}(-\nabla((1-\Delta)^{-1} f^k))$ (see, p.197 of
\cite{kr99}). Hence, we may assume that $f\in \bH^{1}_2(T)$ and show
(\ref{e 6.5.2}) and (\ref{e 6.5.3}) with $\| f\|_{\bH^{1}_2(T)}$ in
place of $\| f\|_{\bH^{-1}_2(T)}$.

 2. By Theorem 4.10 and Theorem 5.1 in \cite{Kr99}, for each $k$ the equation
$$
du^k=\left(D_i(\delta \cdot\delta_{ij}\delta_{kr}
u^r_{x^ix^j}+\bar{f}^{ik})+ f^k\right)dt +g^k_m dw^m_t,
$$
or equivalently,
$$
du^k=(\delta\Delta u^k+\bar{f}^{ik}_{x^i}+f^k)dt+g^k_mdw^m_t, \quad
u^k(0)=u^k_0,
$$
has a solution $u^k$ and we have $u:=(u^1,u^2,\cdots,u^{d_1})^*$ as
the unique solution  of
\begin{eqnarray}
du=(\delta\Delta u +\bar{f}^{i}_{x^i}+ f)dt+g_m dw^m_t, \quad
u(0)=u_0,\nonumber
\end{eqnarray}
in $\mathcal{H}^{1}_{2}(T)$ with  estimates (\ref{e 6.5.2}) and (\ref{e 6.5.3}).  For
$\lambda\in [0,1]$ we define
\begin{eqnarray}
E^{ij}_{\lambda}=(e^{ij}_{kr,\lambda})&:=&(1-\lambda)\left(A^{ij}-\mathcal{A}^{ij}\right)+\lambda
\delta\cdot\delta_{ij}
I\nonumber\\
&=& \left((1-\lambda)A^{ij}+\lambda\delta\cdot\delta_{ij}
I\right)-(1-\lambda)\mathcal{A}^{ij}=A^{ij}_{\lambda}-\mathcal{A}^{ij}_{\lambda},\nonumber
\end{eqnarray}
where $A^{ij}_{\lambda}:=(1-\lambda)A^{ij}+\lambda\delta
\cdot\delta_{ij}I,\;
\mathcal{A}^{ij}_{\lambda}:=(1-\lambda)\mathcal{A}^{ij}$. Then we
have
\begin{eqnarray}
|A^{ij}_{\lambda}|\le |A^{ij}|,\quad |\mathcal{A}^{ij}_{\lambda}|\le
|\mathcal{A}^{ij}|,\quad \delta|\xi|^2\leq
\sum_{i,j}\xi^*_iE^{ij}_{\lambda}\xi_j\nonumber
\end{eqnarray}
for any real $d_1\times d$-matrix $\xi$. Also, we define
$$
\bar{B}^i:=(1-\lambda)\bar{B}^i, \quad
B^i_{\lambda}:=(1-\lambda)B^i, \quad C_{\lambda}:=(1-\lambda)C,
\quad  \mathcal{N}_{\lambda}:=(1-\lambda)\mathcal{N}.
$$
Then $\bar{B}^i_{\lambda},
B^i_{\lambda},C_{\lambda},\mathcal{N}_{\lambda}$ satisfy
(\ref{extra}). Thus, having the method of continuity in mind, we
only prove that  (\ref{e 6.5.2}) and (\ref{e 6.5.3})
hold given that a solution $u$ already exists.

3. Applying the stochastic product rule $d|u^k|^2=2u^kdu^k+du^kdu^k$
for each $k$, we have
\begin{eqnarray}
|u^k(t)|^2&=&\quad |u^k_0|^2\nonumber\\&&+\int^t_0
2u^k\left(D_i(a^{ij}_{kr}u^r_{x^j}+\bar{f}^{ik})+b^i_{kr}u^r_{x^i}+c_{kr}u^r+ f^k\right)ds\nonumber\\
&&+\int^t_0|\sigma^i_{kr}u^r_{x^i}+\nu_{kr}u^r+g^{k}|_{\ell_2}^2ds\nonumber\\
&&+\int^t_0
2u^k(\sigma^i_{kr,m}u^r_{x^i}+\nu_{kr,m}u^r+g^{k}_m)dw^{m}_s,\quad
t>0.\label{square}
\end{eqnarray}
Note that, making the summation on $r,i$ appeared, we get
\begin{eqnarray*}
&&\sum_k\left|\sum_{r,i}\sigma^i_{kr}u^r_{x^i}+\sum_r\nu_{kr}u^r+g^{k}\right|_{\ell_2}^2\\
&=&
2\sum_{i,j}(u_{x^i})^*\mathcal{A}^{ij}u_{x^j}+\sum_k\left[\left|(\mathcal{N}u)^k\right|_{\ell_2}^2
+|g^k|^2_{\ell_2}\right]\\
&&+2\sum_k\left[(\sum_i(\Sigma^iu_{x^i})^k,g^k)_{\ell_2}+((\mathcal{N}u)^k,g^k)_{\ell_2}
+(\sum_i(\Sigma^iu_{x^i})^k,(\mathcal{N}u)^k)_{\ell_2}\right].
\end{eqnarray*}
By taking expectation, integrating with respect to $x$, and using
integrating by parts in turn on (\ref{square}), we obtain
\begin{eqnarray}
&&\bE\int_{\bR^d}|u(t)|^2dx+2\;\bE\int^t_0\int_{\bR^d}\sum_{i,j}(u_{x^i})^*(A^{ij}-\mathcal{A}^{ij})u_{x^j}dxds\nonumber\\
&=&\bE\int_{\bR^d}|u_0|^2dx\nonumber\\
&& +2\sum_i\bE\int^t_0\int_{\bR^d}\left[-2u^*_{x^i}\bar{f}^i+u^*(
B^iu_{x^i})\right]dxds
+2\bE\int^t_0\int_{\bR^d}\left[Cu+u^* f\right]dxds\nonumber\\
&&
+\sum_k\bE\int^t_0\int_{\bR^d}\left[\left|(\mathcal{N}u)^k\right|_{\ell_2}^2
+|g^k|^2_{\ell_2}\right]dxds\nonumber\\
&&+2\sum_k\bE\int^t_0\int_{\bR^d}\left[(\sum_i(\Sigma^iu_{x^i})^k,g^k)_{\ell_2}+((\mathcal{N}u)^k,g^k)_{\ell_2}
+(\sum_i(\Sigma^iu_{x^i})^k,(\mathcal{N}u)^k)_{\ell_2}\right]dxds.\nonumber
\end{eqnarray}
Note that we have
\begin{eqnarray}
2\left|\sum_{k}(\sum_i(\Sigma^iu_{x^i})^k,g^k)_{\ell_2}\right|&\le&
2\sum_{k}\big|\sum_{r,i}\sigma^i_{kr}u^r_{x^i}\big|_{\ell_2}\left|g^k\right|_{\ell_2}\nonumber\\
&\le&\sum_{k}\left(\frac{\varepsilon}{2}
\big|\sum_{r,i}\sigma^i_{kr}u^r_{x^i}\big|_{\ell_2}^2+\frac{2}{\varepsilon}\left|g^k\right|_{\ell_2}^2\right)\nonumber\\
&\le& \frac{\varepsilon}{2}
|u_x|^2\sum_{k,r,i}\big|\sigma^i_{kr}\big|_{\ell_2}^2+\frac{2}{\varepsilon}\sum_k\left|g^k\right|_{\ell_2}^2\nonumber\\
&=& \varepsilon
|u_x|^2\sum_{r,i}\big|\alpha^{ii}_{rr}\big|^2+\frac{2}{\varepsilon}\sum_k\left|g^k\right|_{\ell_2}^2\label{2009.09.25
1}
\end{eqnarray}
for any $\varepsilon>0$; similarly, we get
\begin{eqnarray}
2\left|\sum_{k}((\mathcal{N}u)^k,g^k)_{\ell_2}\right|&\le&
|\mathcal{N}||u|+\sum_k |g^k|^2_{\ell_2},\label{2009.09.25
2}\\
2\left|\sum_{k}(\sum_i(\Sigma^iu_{x^i})^k,(\mathcal{N}u)^k)_{\ell_2}\right|&\le&
\varepsilon
|u_x|^2\sum_{r,i}\big|\alpha^{ii}_{rr}\big|^2+\frac{2}{\varepsilon}|\mathcal{N}||u|.\label{2009.09.25
3}
\end{eqnarray}

Hence, it  follows that
\begin{eqnarray}
&&\bE\int_{\bR^d}|u(t)|^2dx+2\delta\;
\bE\int^t_0\int_{\bR^d}|u_x|^2dxds\nonumber\\
&\leq&\bE\int_{\bR^d}|u_{0}|^2dx +c\varepsilon\;\bE\int^t_0\int_{\bR^d}|u_x|^2dxds+c\bE\int^t_0\int_{\bR^d}|u(s)|^2dxds\nonumber\\
&&+ c\sum_i\bE\int^t_0\int_{\bR^d}|\bar{f}^i|^2dxds
+\bE\int^t_0\int_{\bR^d}|f|^2dxds+c\;\bE\sum_k\int^t_0\int_{\bR^d}|g^{k}|_{\ell_2}^2dxds\nonumber\\
&\le&c\varepsilon\;\bE\int^t_0\int_{\bR^d}|u_x|^2dxds+c\bE\int^t_0\int_{\bR^d}|u(s)|^2dxds\nonumber\\
&&+c\sum_i\|\bar{f}^i\|^2_{\bL_2(T)}+\|
f\|^2_{\bL_2(T)}+c\|g\|^2_{\bL_2(T,\ell_2)}+\|u_0\|^2_{U^{1}_2}.\nonumber
\end{eqnarray}

Choosing small $\varepsilon$, we obtain
\begin{eqnarray}
\|u_x\|^2_{\bL_2(T)}&\le& c(\|u\|^2_{\bL_2(T)}+\|f\|^2_{\bL_2(T)}
+\sum_i\|\bar{f}^i\|^2_{\bL_2(T)}+\|g\|^2_{\bL_2(T,\ell_2)}+\|u_0\|^2_{U^{1}_2}),\nonumber
\\
\bE\int_{\bR^d}|u(t)|^2dx&\le&
c\bE\int^t_0\int_{\bR^d}|u(s)|^2dxds\nonumber\\&&+c( \|
f\|^2_{\bL_2(T)}
+\sum_i\|\bar{f}^i\|^2_{\bL_2(T)}+\|g\|^2_{\bL_2(T,\ell_2)}+\|u_0\|^2_{U^{1}_2}),\nonumber
\end{eqnarray}
where $c$ does not depend on $T$. Now we recall the remark in step 1, and see that the first inequality implies (\ref{e 6.5.2}). Also the second inequality and Gronwall's inequality
lead us to (\ref{e 6.5.3}). The
theorem is proved.
\end{proof}


\mysection{The system on $\bR^d_+$} \label{section half}

In this section we present some results for the systems defined on
$\bR^d_+$. In the next section, these results will be modified and
be used to develop our theory of the systems defined on
$C^1$-domains.

Here we use the Banach spaces  introduced in \cite{kr99}. Let
$\zeta\in C^{\infty}_{0}(\bR_{+})$ be a   function satisfying
\begin{equation}
                                       \label{eqn 5.6.5}
\sum_{n=-\infty}^{\infty}\zeta(e^{n+x})>c>0, \quad \forall x\in \bR,
\end{equation}
where $c$ is a constant. Note that  any nonnegative function
$\zeta$, $\zeta>0$ on $[1,e]$, satisfies (\ref{eqn 5.6.5}). For
$\theta,\gamma \in \bR$, we let $H^{\gamma}_{p,\theta}$ denote the set
of all distributions $u=(u^1,u^2,\cdots u^{d_1})$  on $\bR^d_+$ such
that
\begin{equation}
                                                 \label{10.10.03}
\|u\|_{H^{\gamma}_{p,\theta}}^{p}:= \sum_{n\in\bZ} e^{n\theta}
\|\zeta(\cdot)u(e^{n} \cdot)\|^p_{H^{\gamma}_p} < \infty.
\end{equation}
 If $g=(g^1,g^2,\ldots,g^{d_1})$ and each $g^k$ is an
 $\ell_2$-valued function,
 then we define
$$
\|g\|_{H^{\gamma}_{p,\theta}(\ell_2)}^{p}= \sum_{n\in\bZ}
e^{n\theta} \|\zeta(\cdot)g(e^{n} \cdot)\|^p_{H^{\gamma}_p(\ell_2)}.
$$
It is known (see \cite{kr99}) that up to equivalent norms the space
$H^{\gamma}_{p,\theta}$ is independent of the choice of $\zeta$.
Also,   for any $\eta\in C^{\infty}_0(\bR_+)$, we have
\begin{equation}
                            \label{eqn 5.6.1}
\sum_{n=-\infty}^{\infty}
e^{n\theta}\|u(e^n\cdot)\eta\|^p_{H^{\ga}_p} \leq c
\sum_{n=-\infty}^{\infty}
e^{n\theta}\|u(e^n\cdot)\zeta\|^p_{H^{\ga}_p},
\end{equation}
where $c$ depends only on $d,d_1,\gamma,\theta,p,\eta,\zeta$.
Furthermore,
 if $\gamma$ is a nonnegative integer, then
\begin{equation}
                   \label{eqn compare1}
\|u\|^p_{H^{\gamma}_{p,\theta}}\sim \sum_{n=0}^{\gamma}\sum_{|\alpha|=n}
\int_{\bR^d_+}|(x^1)^n D^{\alpha}u(x)|^p(x^1)^{\theta-d} \,dx.
\end{equation}

Below we collect some other properties of spaces
$H^{\gamma}_{p,\theta}$. Let $M^{\alpha}$ be the operator of
multiplying by $(x^1)^{\alpha}$ and $M=M^1$.

\begin{lemma} $(\cite{kr99})$ Let $d-1<\theta<d-1+p$.
              \label{lemma collection}

 (i) Assume that $\gamma-d/p=m+\nu$ for some
$m=0,1,\cdots$ and $\nu\in (0,1]$.  Then for any $u\in
H^{\gamma}_{p,\theta}$ and $i\in \{0,1,\cdots,m\}$, we have
$$
|M^{i+\theta/p}D^iu|_{C}+[M^{m+\nu+\theta/p}D^m u]_{C^{\nu}}\leq c
\|u\|_{ H^{\gamma}_{p,\theta}}.
$$

(ii) Let $\alpha\in \bR$, then
$M^{\alpha}H^{\gamma}_{p,\theta+\alpha p}=H^{\gamma}_{p,\theta}$,
$$
\|u\|_{H^{\gamma}_{p,\theta}}\leq c
\|M^{-\alpha}u\|_{H^{\gamma}_{p,\theta+\alpha p}}\leq
c\|u\|_{H^{\gamma}_{p,\theta}}.
$$

(iii) $M D, DM: H^{\gamma}_{p,\theta}\to H^{\gamma-1}_{p,\theta}$
are bounded linear operators.

(iv) There is a constant $c=c(d,p,\theta,\gamma)>0$ so that
$$
c^{-1}\|M^{-1}u\|_{H^{\gamma}_{p,\theta}}\leq
\|u_x\|_{H^{\gamma-1}_{p,\theta}}\leq c
\|M^{-1}u\|_{H^{\gamma}_{p,\theta}}.
$$

\end{lemma}

We define the following stochastic Banach spaces.
$$
\bH^{\gamma}_{p,\theta}(T)=L_p(\Omega\times[0,T],\mathcal{P},H^{\gamma}_{p,\theta}),
\quad
\bH^{\gamma}_{p,\theta}(T,\ell_2)=L_p(\Omega\times[0,T],\mathcal{P},H^{\gamma}_{p,\theta}(\ell_2))
$$
$$
\bL_{p,\theta}(T):=\bH^{0}_{p,\theta}(T),\quad
\bL_{p,\theta}(T,\ell_2):=\bH^{0}_{p,\theta}(T,\ell_2),\quad
U^{\gamma}_{p,\theta}=L_p(\Omega,\mathcal{F}_0,M^{1-2/p}H^{\gamma-2/p}_{p,\theta}).
$$

\begin{definition}
We write $u\in \frH^{\gamma+2}_{p,\theta}(T)$ if $u\in
M\bH^{\gamma+2}_{p,\theta}(T)$, $u(0)\in U^{\gamma+2}_{p,\theta}$
and for some $f\in M^{-1}\bH^{\gamma}_{p,\theta}(T)$, $g\in
\bH^{\gamma+1}_{p,\theta}(T,\ell_2)$,
$$
du=fdt+g_mdw^m_t
$$
holds in the sense of the distributions. The norm in
$\frH^{\gamma+2}_{p,\theta}(T)$ is defined by
$$
\|u\|_{\frH^{\gamma+2}_{p,\theta}(T)}=\|M^{-1}u\|_{\bH^{\gamma+2}_{p,\theta}(T)}+\|Mf\|_{\bH^{\gamma}_{p,\theta}(T)}+
\|g\|_{\bH^{\gamma+1}_{p,\theta}(T,\ell_2)}
+\|u(0)\|_{U^{\gamma+2}_{p,\theta}}.
$$
\end{definition}

Let us denote
 $$
 K:=\sqrt{\sum_j (K^j)^2}.
 $$
\begin{lemma}
                \label{a priori 1}
Assume
\begin{equation}
                     \label{theta 1}
\theta\in \left(d-\frac{\delta}{2K-\delta},\,\,
d+\frac{\delta}{2K+\delta}\right),
\end{equation}
 $\bar{b}^i=b^i=c=0$ and $\nu=0$. Then if    $u\in
M\bH^{1}_{2,\theta}(T)$ is a solution of  system (\ref{eqn main
system}) on $[0,T]\times \mathbb{R}^d_+$ and $u\in
L_2(\Omega,C([0,T],C^1_0((1/N,N)\times \{x':|x'|<N\})))$ for some
$N>0$, then we have
\begin{equation}
               \label{eqn main}
\|M^{-1}u\|^2_{\bH^1_{2,\theta}(T)}\leq
c(\|\bar{f}^i\|_{\bL_{2,\theta}(T)}+\|M
f\|^2_{\bH^{-1}_{2,\theta}(T)}+\|g\|^2_{\bL_{2,\theta}(T,\ell_2)}+\|u_0\|^2_{U^1_{2,\theta}}),
\end{equation}
where $c=c(d,d_1,\delta,\theta,K,L)$.
\end{lemma}

\begin{proof}
1. By Corollary 2.12 in \cite{kr99}, $f^k$ has the following
representation:
$$
f^k=\sum_{i=1}^d D_iF^{ik}, \quad
\sum_i\|F^{ik}\|_{L_{2,\theta}}\leq c\|Mf^k\|_{H^{-1}_{2,\theta}}.
$$
Also since $ \|M^{-1}u\|_{H^1_{2,\theta}}\leq
c\|u_x\|_{L_{2,\theta}}$ (\;see Lemma \ref{lemma collection}(iv)\;),
 it is enough to assume
$f^k=0$ and  prove
$$
\|u_x\|^2_{\bL_{2,\theta}(T)}\leq
c(\|\bar{f}^i\|_{\bL_{2,\theta}(T)}+\|g\|^2_{\bL_{2,\theta}(T,\ell_2)}+\|u_0\|^2_{U^1_{2,\theta}}).
$$
2. Again, as in the proof of Theorem \ref{thm 1}, applying the
stochastic product rule  $d|u^k|^2=2u^kdu^k+du^kdu^k$ for each $k$,
we get
\begin{eqnarray}
|u^k(t)|^2&=&|u^k_0|^2+\int^t_0
2u^k \left[D_i(a^{ij}_{kr}u^r_{x^j}+\bar{f}^{ik})\right]ds\nonumber\\
&&+\int^t_0|\sigma^i_{kr}u^r_{x^i}+g^{k}|_{\ell_2}^2ds +\int^t_0
2u^k(\sigma^i_{kr,m}u^r_{x^i}+g^{k}_m)dw^{m}_s,\nonumber
\end{eqnarray}
where the summations on $i,j,r$ are understood. Denote $c=\theta-d$.
For each  $k$, we have
\begin{eqnarray}
0&\leq& \bE\int_{\bR^d_+}|u^k(T,x)|^2(x^1)^c
dx\nonumber\\
&=&\bE\int_{\bR^d_+}|u^k(0,x)|^2(x^1)^cdx\nonumber\\
&& + 2\bE\int^T_0\int_{\bR^d_+} u^kD_i(a^{ij}_{kr}u^r_{x^j})(x^1)^c
dxds+2\bE\int^T_0\int_{\bR^d_+} u^r\bar{f}^{ik}_{x^i}(x^1)^c
dxds\nonumber\\
&&+
 \bE\int^T_0\int_{\bR^d_+}
|\sigma^i_{kr}u^r_{x^i}|^2_{\ell_2}(x^1)^c dxds\nonumber\\
&&+2\bE\int^T_0\int_{\bR^d_+}(\Sigma^iu_{x^i})^k,g^k)_{\ell_2}
(x^1)^c dxds + \bE\int^T_0\int_{\bR^d_+} |g^k|^2_{\ell_2}(x^1)^c
dxds.\label{2009.06.02 04:27 PM}
\end{eqnarray}
Note that, by integration by parts, we get
$$
2\bE\int^T_0\int_{\bR^d_+} u^r\bar{f}^{ik}_{x^i}(x^1)^c dxds
=-2\bE\int^T_0\int_{\bR^d_+} \left[u^r_{x^i}\bar{f}^{ik}(x^1)^c
+cM^{-1}u^r\bar{f}^{1k}(x^1)^c\right]dxds
$$
$$
\leq \varepsilon \|u_x\|^2_{\bL_{2,\theta}(T)}+
\varepsilon\|M^{-1}u\|^2_{\bL_{2,\theta}(T)}+c(\varepsilon)\|\bar{f}\|^2_{\bL_{2,\theta}(T)}.
$$
Also, the second term in the right hand side of (\ref{2009.06.02
04:27 PM}) is
$$
\bE\int^T_0\int_{\bR^d_+}\left[-2a^{ij}_{kr}u^k_{x^i}u^r_{x^j}-2
c(a^{1j}_{kr}u^r_{x^j})(M^{-1}u^k)\right](x^1)^cdxds.
$$
 Thus, by summing up the terms in
(\ref{2009.06.02 04:27 PM}) over $k$ and rearranging the terms, we
obtain
\begin{eqnarray}
&&2\bE\int^T_0\int_{\bR^d_+}u^*_{x^i}\left(A^{ij}-\mathcal{A}^{ij}\right)u_{x^j}\;(x^1)^c
dxds\nonumber\\
&\leq&|c|\left(\kappa\|u_x\|^2_{\bL_{2,\theta}(T)}+
K^2\kappa^{-1}\|M^{-1}u\|^2_{\bL_{2,\theta}(T)}\right)\nonumber\\
&&+N\varepsilon\left(\|M^{-1}u\|^2_{\bL_{2,\theta}(T)}
+\|u_x\|^2_{\bL_{2,\theta}(T)}\right)\nonumber\\
&&+c(\varepsilon)\left(\|\bar f^i\|^2_{\bL_{2,\theta}(T)}
+\|g\|^2_{\bL_{2,\theta}(T,\ell_2)}\right)+\|u(0)\|^2_{U^1_{2,\theta}},\label{eqn
2}
\end{eqnarray}
for any
$\kappa,\varepsilon>0$. This is because for any vectors $v,w\in \bR^n$ and
$\kappa>0$
$$
|<A^{1j}v,w>|\leq |A^{1j}v||w|\leq K^j|v||w|\leq
\frac12(\kappa|v|^2+\kappa^{-1}(K^j)^2|w|^2)
$$
and consequently,
\begin{eqnarray}
&&\bE\int^T_0\int_{\bR^d_+}\left[-2a^{ij}_{kr}u^k_{x^i}u^r_{x^j}-2
c(a^{1j}_{kr}u^r_{x^j})(M^{-1}u^k)\right](x^1)^cdxds \nonumber\\
&\leq& \bE\int^T_0\int_{\bR^d_+}-2a^{ij}_{kr}u^k_{x^i}u^r_{x^j}\,dxds+
|c|\left(\kappa\|u_x\|^2_{\bL_{2,\theta}(T)}+
K^2\kappa^{-1}\|M^{-1}u\|^2_{\bL_{2,\theta}(T)}\right). \label{2009.09.02.6:18PM}
\end{eqnarray}
Now, Assumption
(\ref{assumption 1}),  inequality (\ref{eqn 2}), the inequality
\begin{equation}
                    \label{eqn 3}
\|M^{-1}u\|^2_{L_{2,\theta}}\leq
\frac{4}{(d+1-\theta)^2}\|u_x\|^2_{L_{2,\theta}}
\end{equation}
(see Corollary 6.2 in \cite{kr99}), and Lemma \ref{lemma collection} (iv)
lead us to
$$
 2\delta\|u_x\|^2_{\bL_{2,\theta}(T)}-
 |c|\left(\kappa
 +\frac{4K^2}{\kappa(d+1-\theta)^2}\right)\|u_x\|^2_{\bL_{2,\theta}(T)}
 $$
 $$
\leq N\varepsilon \|u_x\|^2_{\bL_{2,\theta}(T)}
+N\|\bar{f}^i\|^2_{\bL_{2,\theta}(T)}+N\|g\|^2_{\bL_{2,\theta}(T)}+\|u(0)\|^2_{U^1_{2,\theta}}.
$$
Now, it is enough to take $\kappa=2K/(d+1-\theta)$ and observe that
(\ref{theta 1}) is equivalent to the condition
$$
2\delta-|c|\left(\kappa
 +\frac{4K}{\kappa(d+1-\theta)^2}\right)=2\delta-\frac{4|c|K}{d+1-\theta}>0.
 $$
The lemma is proved.
\end{proof}

Here is the main result of this section.

\begin{theorem}
                 \label{theorem half}
Suppose (\ref{theta 1}) holds and
\begin{equation}
                 \label{eqn 9.16.1}
 |M\bar{b}^i|+|Mb^i|+|M^2c|+|M\nu|_{\ell_2}<\beta.
\end{equation}
 Then there
exists constant $\beta_0=\beta_0(d,d_1,\theta,\delta,  K,L)>0$ so that if
$\beta\leq \beta_0$, then for any $\bar{f}^i\in \bL_{2,\theta}(T)$,
$f\in M^{-1}\bH^{-1}_{2,\theta}(T)$, $g\in
\bL_{2,\theta}(T,\ell_2)$, and $u_0\in U^{1}_{2,\theta}$, system
(\ref{eqn main system}) has a unique solution $u\in
\frH^{1}_{2,\theta}(T)$, and furthermore
\begin{equation}
                                  \label{eqn 9.9.2}
 \|u\|_{\frH^{1}_{2,\theta}(T)}\leq
 c\|\bar{f}^i\|_{\bL_{2,\theta}(T)}+
 c\|Mf\|_{\bH^{-1}_{2,\theta}(T)}
 +c\|g\|_{\bL_{2,\theta}(T,\ell_2)}+
 c\|u_0\|_{U^{1}_{2,\theta}}
 \end{equation}
 where $c=c(d,\delta,K,L,T)$.
 \end{theorem}

\begin{proof}
 As before, we only prove that the a priori estimate (\ref{eqn 9.9.2}) holds given that
 a solution $u$ already exists. By Theorem 2.9 in \cite{KL2}, for any nonnegative integer $n\geq
\gamma$, the set
$$
\frH^{n}_{2,\theta}(T) \cap
\bigcup_{N=1}^{\infty}L_2(\Omega,C([0,T],C^n_0((1/N,N)\times
\{x':|x'|<N\})))
$$
is everywhere dense in $\frH^{\gamma}_{2,\theta}(T)$ and thus we may
assume that $u$ is sufficiently smooth in $x$ and vanishes near the
boundary.

{\bf{Step 1}}. If $\bar{b}^i=b^i=c=0$ and $\nu=0$ the  a priori
estimate follows from Lemma \ref{a priori 1}.

{\bf{Step 2}}. In general, by Step 1,
\begin{eqnarray*}
\|M^{-1}u\|_{\bH^{1}_{2,\theta}(T)} & \leq &
 c\|M\bar{b}^iM^{-1}u+\bar{f}^i\|_{\bL_{2,\theta}(T)}+
 c\|Mb^iu_{x^i}+M^2cM^{-1}u+Mf\|_{\bH^{-1}_{2,\theta}(T)}\\
 &&
 \quad +c\|M\nu M^{-1}u+g\|_{\bL_{2,\theta}(T,\ell_2)}+
 c\|u_0\|_{U^{1}_{2,\theta}}.
 \end{eqnarray*}
 Since $\|\cdot\|_{H^{-1}_{2,\theta}}\leq \|\cdot
 \|_{L_{2,\theta}}$,
 we easily see that the above is less than
 $$
 c\beta \|M^{-1}u\|_{\bH^1_{2,\theta}(T)}+c\|\bar{f}\|_{\bL_{2,\theta}(T)}+
 c\|Mf\|_{\bH^{-1}_{2,\theta}(T)}
 +c\|g\|_{\bL_{2,\theta}(T,\ell_2)}+
 c\|u_0\|_{U^{1}_{2,\theta}}.
 $$
 Now it is enough to take $\beta_0$ so that $c \beta <1/2$ for any
 $\beta\leq \beta_0$. The theorem is proved.

\end{proof}

\begin{remark}
We do not know how sharp  (\ref{theta 1}) is. However, if $\theta \not\in (d-1,d+1)$ then Theorem \ref{theorem half} is false even for the heat equation $u_t=\Delta u+f$  (see \cite{kr99}).

We also mention that if the coefficients are sufficiently smooth in $x$,
then one can get quite wider range of $\theta$. This will be shown in the subsequent article \cite{KL4}.

\end{remark}

\mysection{The system on  $\mathcal{O} \subset \bR^d$}
\label{section domain}

In this section we assume the following.

\begin{assumption}
                                         \label{assumption domain}

The domain $\cO$  is of class $C^{1}_{u}$. In other words, for any
$x_0 \in \partial \cO$, there exist constants $r_0, K_0\in(0,\infty)
$ and  a one-to-one continuously differentiable mapping $\Psi$ of
 $B_{r_0}(x_0)$ onto a domain $J\subset\bR^d$ such that

(i) $J_+:=\Psi(B_{r_0}(x_0) \cap \cO) \subset \bR^d_+$ and
$\Psi(x_0)=0$;

(ii)  $\Psi(B_{r_0}(x_0) \cap \partial \cO)= J \cap \{y\in
\bR^d:y^1=0 \}$;

(iii) $\|\Psi\|_{C^{1}(B_{r_0}(x_0))}  \leq K_0 $ and
$|\Psi^{-1}(y_1)-\Psi^{-1}(y_2)| \leq K_0 |y_1 -y_2|$ for any $y_i
\in J$;

(iv)   $\Psi_{x}$ is uniformly continuous in for $B_{r_{0}}(x_{0})$.
\end{assumption}

To proceed further we introduce some well known results from
\cite{GH} and \cite{KK2}.

\begin{lemma}
                                           \label{lemma 10.3.1}
Let the domain $\cO$ be of class $C^{1}_{u}$. Then

(i) there is a bounded real-valued function $\psi$ defined in
$\bar{\cO} $  such that the functions $\psi(x)$ and
$\rho(x):=\text{dist}(x,\partial \cO)$ are comparable  in the part
of
  a neighborhood of $\partial \cO$ lying in $\cO$. In other words, if $\rho(x)$ is
sufficiently small, say $\rho(x)\leq 1$, then $N^{-1}\rho(x) \leq
\psi(x) \leq N\rho(x)$ with some constant
 $N$ independent of $x$,

 (ii) for any  multi-index $\alpha$ it holds that
\begin{equation}
                                                             \label{03.04.01}
\sup_{\cO} \psi ^{|\alpha|}(x)|D^{\alpha}\psi_{x}(x)| <\infty.
\end{equation}

\end{lemma}

To describe the assumptions of $\bar{f}^i$s, $f$, and $g$ in
(\ref{eqn main system}) with space domain $\cO$ we use the Banach
spaces introduced in \cite{KK2} and \cite{Lo2}.
 Let  $\zeta\in C^{\infty}_{0}(\bR_{+})$
be a   nonnegative function satisfying (\ref{eqn 5.6.5}). For $x\in
\cO$ and $n\in\bZ:=\{0,\pm1,...\}$ we define
$$
\zeta_{n}(x)=\zeta(e^{n}\psi(x)).
$$
Then  we have $\sum_{n}\zeta_{n}\geq c$ in $\cO$ and
\begin{equation*}
\zeta_n \in C^{\infty}_0(\cO), \quad |D^m \zeta_n(x)|\leq
N(m)e^{mn}.
\end{equation*}

For $\theta,\gamma \in \bR$, let $H^{\gamma}_{p,\theta}(\cO)$ denote
the set of all distributions $u=(u^1,u^2,\cdots u^{d_1})$  on $\cO$
such that
\begin{equation}
                                                 \label{10.10.03}
\|u\|_{H^{\gamma}_{p,\theta}(\cO)}^{p}:= \sum_{n\in\bZ} e^{n\theta}
\|\zeta_{-n}(e^{n} \cdot)u(e^{n} \cdot)\|^p_{H^{\gamma}_p} < \infty.
\end{equation}
 If $g=(g^1,g^2,\ldots,g^{d_1})$ and each $g^k$ is an
 $\ell_2$-valued function,
 then we define
$$
\|g\|_{H^{\gamma}_{p,\theta}(\cO,\ell_2)}^{p}= \sum_{n\in\bZ}
e^{n\theta} \|\zeta_{-n}(e^{n} \cdot)g(e^{n}
\cdot)\|^p_{H^{\gamma}_p(\ell_2)}.
$$
It is known (see, for instance, \cite{Lo2}) that up to equivalent
norms the space $H^{\gamma}_{p,\theta}(\cO)$ is independent of the
choice of $\zeta$ and $\psi$. Moreover, if $\gamma=n$ is a
non-negative integer, then it holds that
\begin{equation}
                        \label{eqn compare2}
\|u\|^p_{H^{\gamma}_{p,\theta}(\cO)} \sim \sum_{k=0}^n
\sum_{|\alpha|=k}\int_{\cO} |\psi^kD^{\alpha}u(x)|^p
\psi^{\theta-d}(x) \,dx.
\end{equation}
By comparing (\ref{eqn compare1}) and (\ref{eqn compare2}), one
finds that   two spaces $H^{\gamma}_{p,\theta}(\bR^d_+)$ and
$H^{\gamma}_{p,\theta}$ are different since $\psi$ is bounded. Also,
it is easy to see that, for any nonnegative function
$\xi=\xi(x^1)\in C^{\infty}_0(\bR^1)$ satisfying  $\xi=1$ near
$x^1=0$, we have
\begin{equation}
                             \label{eqn 9.9}
\|u\|_{H^{\gamma}_{p,\theta}(\bR^d_+)}\sim \left(\|\xi
u\|_{H^{\gamma}_{p,\theta}}+\|(1-\xi)u\|_{H^{\gamma}_p}\right).
\end{equation}
In particular, if $u(x)=0$ for $x\geq r$, then for any $\alpha\in
\bR$ we get
\begin{equation}
                      \label{eqn 9.99}
c^{-1}\|M^{\alpha}u\|_{H^{\gamma}_{p,\theta}}\leq \|\psi^{\alpha}
u\|_{H^{\gamma}_{p,\theta}(\bR^d_+)}\leq
c\|M^{\alpha}u\|_{H^{\gamma}_{p,\theta}},
\end{equation}
where $c=c(r,\alpha,\gamma,p,\theta)$. We also mention that the
space $H^{\gamma}_{p,\theta}$ can be defined on the basis of
(\ref{10.10.03}) by  formally taking $\psi(x)=x^{1}$ so that
$\zeta_{-n}(e^{n} x)=\zeta(x)$ and (\ref{10.10.03}) becomes
\begin{equation*}
\|u\|^p_{H^{\gamma}_{p,\theta}}:= \sum_{n\in\bZ} e^{n\theta}
\|u(e^{n} \cdot)\zeta \|^p_{H^{\gamma}_p} < \infty.
\end{equation*}

We place the following lemma similar to Lemma \ref{lemma
collection}.
\begin{lemma} $(\cite{kr99})$ Let $d-1<\theta<d-1+p$.

                  \label{lemma collection2}
 Assertions (i)-(iii) in Lemma \ref{lemma collection} hold true with $\psi$
and $H^{\gamma}_{p,\theta}(\cO)$ in place of $M$ and
$H^{\gamma}_{p,\theta}$, respectively.

\end{lemma}

We define
$$
\bH^{\gamma}_{p,\theta}(\cO,T)=L_p(\Omega\times
[0,T],\cP,H^{\gamma}_{p,\theta}(\cO)), \quad
\bH^{\gamma}_{p,\theta}(\cO,T,\ell_2)=L_p(\Omega\times
[0,T],\cP,H^{\gamma}_{p,\theta}(\cO,\ell_2)),
$$
$$
U^{\gamma}_{p,\theta}(\cO)=
\psi^{1-2/p}L_{p}(\Omega,\cF_0,H^{\gamma-2/p}_{p,\theta}(\cO)),
\quad \bL_{p,\theta}(\cO,T)=\bH^{0}_{p,\theta}(\cO,T).
$$
\begin{definition}
We define $ \frH^{\gamma+2}_{p,\theta}(\cO,T)$ as the space of all
functions
 $u=(u^1,\cdots, u^{d_1})\in \psi\bH^{\gamma+2}_{p,\theta}(\cO,T)$ such that
$u(0,\cdot) \in U^{\gamma+2}_{p,\theta}(\cO)$ and  for some $f \in
\psi^{-1}\bH^{\gamma}_{p,\theta}(\cO,T)$, $g\in
\bH^{\gamma+1}_{p,\theta}(\cO,T,\ell_2)$,
$$
du= f \,dt + g_m \, dw^m_t,
$$
in the sense of distributions. The norm in  $
\frH^{\gamma+2}_{p,\theta}(\cO,T)$ is introduced by
$$
\|u\|_{\frH^{\gamma+2}_{p,\theta}(\cO,T)}=
\|\psi^{-1}u\|_{\bH^{\gamma+2}_{p,\theta}(\cO,T)} + \|\psi
f\|_{\bH^{\gamma}_{p,\theta}(\cO,T)}  +
\|g\|_{\bH^{\gamma+1}_{p,\theta}(\cO,T,\ell_2)} +
\|u(0,\cdot)\|_{U^{\gamma+2}_{p,\theta}(\cO)}.
$$

\end{definition}

The following result  is due to N.V.Krylov (see, for instance, \cite{Kr01}).

\begin{lemma}
                             \label{lemma 15.05}
 Let $p\geq 2$. Then there exists a constant
$c=c(d,p,\theta,\gamma,T)$ such that
$$
\bE \sup_{t\leq T}\|u(t)\|^p_{H^{\gamma+1}_{p,\theta}(\cO)}\leq c
\|u\|^p_{\frH^{\gamma+2}_{p,\theta}(\cO,T)}.
$$
In particular, for any $t\leq T$,
$$
\|u\|^p_{\bH^{\gamma+1}_{p,\theta}(\cO,t)}\leq c \int^t_0
\|u\|^p_{\frH^{\gamma+2}_{p,\theta}(\cO,s)}ds.
$$
\end{lemma}

\begin{assumption}
             \label{assumption regularity}
 There is   control on the behavior of $\bar{b}^i_{kr}$,
$b^i_{kr}$, $c_{kr}$
 and $\nu_{kr}$ near
$\partial \cO$, namely,
\begin{equation}
                                                       \label{05.04.01}
\lim_{\substack{\rho(x)\to0\\
x\in \cO}} \sup_{t,
\omega}[\rho(x)|\bar{b}^i_{kr}(t,x)|+\rho(x)|b^i_{kr}(t,x)|+\rho^{2}(x)|c_{kr}(t,x)|+
\rho(x)|\nu_{kr}(t,x)|_{\ell_2}]=0.
\end{equation}
\end{assumption}

Note that Assumption \ref{assumption regularity} allows the coefficients to be unbounded and to blow up
near the boundary. (\ref{05.04.01}) holds if, for instance,
$$
|\bar{b}^i_{kr}(t,x)|+|b^i_{kr}(x)|+|\nu_{kr}(x)|_{\ell_2}\leq c\rho^{-1+\varepsilon}(x), \quad |c_{kr}(t,x)|\leq \rho^{-2+\varepsilon}(x),
$$
for some $c,\varepsilon>0$.

Here is the main result of this section.

\begin{theorem}
                    \label{main theorem on domain}
Let $\cO=\bR^d_+$ or $\cO$ be bounded. Suppose (\ref{theta 1})
and Assumption \ref{assumption regularity} hold. Then for any
$\bar{f}^i\in \bL_{2,\theta}(\cO,T)$ (i=1,\ldots,d), $f\in
\psi^{-1}\bH^{-1}_{2,\theta}(\cO,T),\;g\in
\bL_{2,\theta}(\cO,T,\ell_2)$, and $u_0\in U^{1}_{2,\theta}(\cO)$,
the system (\ref{eqn main system}) admits a unique solution $u\in
\frH^{1}_{2,\theta}(\cO,T)$, and for this solution we have
\begin{equation}
                        \label{main a priori}
\|\psi^{-1}u\|_{\bH^{1}_{2,\theta}(\cO,T)}\leq
c\|\bar{f}^i\|_{\bL_{2,\theta}(\cO,T)}+c\|\psi
f\|_{\bH^{-1}_{2,\theta}(\cO,T)}
+c\|g\|_{\bL_{2,\theta}(\cO,T,\ell_2)}
+c\|u_0\|_{U^{1}_{2,\theta}(\cO)},
\end{equation}
where $c=c(d,\delta,\theta,K,L)$.
\end{theorem}

\begin{remark}
By carefully inspecting  our arguments below one can check that
Theorem \ref{main theorem on domain} holds even if the $C^1$-domain
$\cO$ is not bounded.
\end{remark}

To prove Theorem \ref{main theorem on domain} we need the following
a priori estimate near the boundary.

\begin{lemma}
                            \label{lemma a priori}
 Suppose that $u\in \frH^1_{2,\theta}(\cO,T)$ is a solution
of system (\ref{eqn main system}) such that $u(t,x)=0$ for $x\in \cO
\backslash B_r(x_0)$, $x_0\in
\partial \cO$. Then there exists constant
$r_1\in (0,1)$, independent of $x_0$ and $u$, such that if
 $r\leq r_1$, then a priori estimate (\ref{main a priori}) holds.
\end{lemma}

\begin{proof}

 Let $x_0 \in \partial \cO$ and $\Psi$ be a function from
Assumption \ref{assumption domain}. In \cite{KK2} it is shown that
 $\Psi$ can be chosen in such a way that
\begin{equation}
                                                     \label{2.25.02}
\rho(x)\Psi_{xx}(x) \to 0 \quad \text{as}\quad  x\in
B_{r_0}(x_0)\cap \cO,
 \text{and} \,\,\,  \rho(x) \to 0,
\end{equation}
where  the
 convergence in (\ref{2.25.02}) is independent of  $x_0$.

 Define $r=r_{0}/K_{0}$
and fix smooth functions $\eta \in C^{\infty}_{0}(B_r ), \varphi\in
C^{\infty}(\bR)$ such that $ 0 \leq \eta, \varphi \leq 1$, and
 $\eta=1$ in $B_{r/2} $,  $\varphi(t)=1$ for $t\leq -3$, and
$\varphi(t)=0$
 for $t\geq-1$ and  $0\geq\varphi'\geq-1$. We observe that
$\Psi(B_{r_0}(x_0))$ contains $B_r $.
 For $n=1,2,... $,\; $t>0$, $x\in\bR^{d}_{+}$ we introduce
$\varphi_{n}(x):=\varphi(n^{-1}\ln x^1)$,
\begin{eqnarray*}
\hat{a}^{ij}(t,x)&:=& \eta(x)\left(\sum_{l,m=1}^d
  a^{lm}(t,\Psi^{-1}(x))\cdot\partial_l\Psi^{i}(\Psi^{-1}(x))\cdot\partial_m\Psi^{j}(\Psi^{-1}(x))\right)  +
\delta^{ij}(1- \eta(x) )I,\\
\hat{\bar{b}}^{i,n}(t,x)&:=&\eta(x)
\varphi_{n}(x)\sum_{l}\bar{b}^{l}(t,\Psi^{-1}(x))\cdot\partial_l\Psi^{i}(\Psi^{-1}(x)),\\
\hat{b}^{i,n}(t,x) &:=&\eta(x) \varphi_{n}(x)\Big[
-\sum_{l,m,r,j}a^{lm}(t,\Psi^{-1}(x))\cdot
(\partial_m\Psi^{j}\cdot\partial_{lr}\Psi^{i})(\Psi^{-1}(x))\cdot\partial_j(\Psi^{-1})^r(x)\\
&&\qquad\qquad\qquad\qquad+\sum_{l}b^{l}(t,\Psi^{-1}(x))\cdot\partial_l\Psi^{i}(\Psi^{-1}(x))\Big],\\
\hat{c}^{n}(t,x) &:=&\eta(x) \varphi_{n}(x)c(t,\Psi^{-1}(x)),\\
\hat{\sigma}^{i}(t,x)&:=&\eta(x)\sum_{l}\sigma^{l}(t,\Psi^{-1}(x))\cdot\partial_l\Psi^{i}(\Psi^{-1}(x)),\\
\hat{\nu}^n(t,x)&:=&\eta(x) \varphi_{n}(x)\nu(t,x)(t,\Psi^{-1}(x)).
\end{eqnarray*}
Then $(\hat{a}^{ij},\hat{\sigma}^i)$ satisfies (\ref{assumption 1})
and (\ref{assumption 2}). We take  $\beta_0$ from Theorem
\ref{theorem half} corresponding to $d,d_1,\theta,\delta, L$ and
$K$. We observe that $\varphi_n(x)=0$ for $x^1 \geq e^{-n}$. Also,
note
 that
 (\ref{2.25.02}) implies $x^{1}\Psi_{xx}(\Psi^{-1}(x)) \to 0$ as $x^1\to 0$.
  Using these facts and (\ref{05.04.01}),
 one can fix
 $n>0$ which is sufficiently large, independent of $x_0$, and
$$
x^1 |\hat{\bar{b}}^{i,n}_{kr}(t,x)| +x^1 |\hat{b}^{i,n}_{kr}(t,x)| +
(x^1)^2|\hat{c}^{n}_{kr}(t,x)|+ x^1|\hat{\nu}^n_{kr}(t,x)|_{\ell_2} \leq
\beta_0, \quad \forall\;\; \omega,t,x.
$$
  Now, we fix   $r_1  <r_0  $ so that
\begin{equation}
                    \label{eqn 8.21.4}
\Psi(B_{r_1}(x_0)) \subset B_{r/2} \cap \{x:x^1 \leq e^{-3n}\}.
\end{equation}

Next, we observe  that, by Lemma \ref{lemma 10.3.1} and
  Theorem 3.2 in \cite{Lo2} (or see \cite{KK2}),
for any $\nu,\alpha \in \bR $ and $h \in
\psi^{-\alpha}H^{\nu}_{p,\theta}(\cO)$ with support in
$B_{r_0}(x_0)$ we have
\begin{equation}
                                                          \label{1.28.01}
\|\psi^{\alpha}h\|_{H^{\nu}_{p,\theta}(\cO)} \sim
\|M^{\alpha}h(\Psi^{-1})\|_{H^{\nu}_{p,\theta}}.
\end{equation}
Thus, for $v(t,x):=u(t,\Psi^{-1}(x))$ we have $u\in
\frH^1_{2,\theta}(T)$ and $v$ satisfies
\begin{eqnarray*}
dv^k&=&(D_i(\hat{a}^{ij}_{kr}v^r_{x^j}+\hat{\bar{b}}^{i,n}_{kr}v^r+\hat{\bar{f}}^{ik})
+\hat{b}^{i,n}_{kr}v^r_{x^i}+\hat{c}^n_{kr}v^r+ \hat{f}^k)dt\\
&&+(\hat{\sigma}^i_{kr,m}v^r_{x^i}+\hat{\nu}^n_{kr,m}v^r+\hat{g}^k_m)dw^m_t,
\end{eqnarray*}
where
$$
\hat{\bar{f}}^{ik}=\sum_{\ell}(\bar{f}^{ik}\partial_i\Psi^{\ell})(\Psi^{-1}(x)),
\quad
\hat{f}^k=\bar{f}^{ik}(\Psi^{-1}(x))\partial_{ij}\Psi^{\ell}(\Psi^{-1}(x))\partial_i(\Psi^{-1})^j(x)+f^k(\Psi^{-1}(x)).
$$
Hence, the a priori estimate follows from Theorem \ref{theorem half}
and (\ref{1.28.01}). The lemma is proved.

\end{proof}

\begin{remark}
                  \label{remark 10.31}
Let $\cO=\bR^d_+$. Then, in fact,    Lemma \ref{lemma a priori}
holds if $u(t,x)=0$ for $x^1 \geq r_1$ for some $r_1$. Indeed, by
(\ref{05.04.01}) there is $r_1>0$ so that
\begin{equation}
                        \label{eqn 10.31.1}
 |M\bar{b}^i|+|Mb^i|+|M^2c|+|M\nu|_{\ell_2}<\beta_0
\end{equation}
for $x^1\leq r_1$. Now, if $u(t,x)=0$ for $x^1 \geq r_1$, then
without affecting the system we may put $\bar{b}^i=b^i=c=0$ and
$\nu=0$ for $x^1\geq r_1$ so that (\ref{eqn 10.31.1}) holds for all
$x$. Consequently the assertion follows from Theorem \ref{theorem
half} and (\ref{eqn 9.99}).
\end{remark}

Next, we prove the a priori estimate for small $T$.

\begin{lemma}
                \label{short time}
 Let assumptions in Theorem \ref{main theorem on domain} be
 satisfied. Then
 there exists a constant $\varepsilon \in (0,1)$ so that if $T\leq
 \varepsilon$, then a priori estimate (\ref{main a priori}) holds for any solution $u\in \frH^1_{2,\theta}(\cO,T)$
 of system (\ref{eqn main system}) with $u_0=0$.
\end{lemma}

\begin{proof}
We prove the lemma only when $\cO$ is bounded. The case
$\cO=\bR^d_+$ is treated similarly. Take a partition of unity
$\{\zeta_n:n=0,1,2,...,N_0\}$, where $N_0<\infty$, such that $\zeta_0\in C^{\infty}_0(\cO)$
and $\zeta_n\in C^{\infty}_0(B_{r_1/2}(x_n))$ with $x_n\in \partial
\cO$ for $n=1,\ldots,N_0$. Also, we fix functions $\bar{\zeta}_n$ such
that $\bar{\zeta}_0\in C^{\infty}_0(\cO), \bar{\zeta}_n\in
C^{\infty}_0(B_{r_1}(x_n))$ for $n=1,\ldots,N_0$, and
$\zeta_n\bar{\zeta}_n=\zeta_n$ for each $n$. We note that
$v_n:=u\zeta_n$ satisfies
\begin{eqnarray}
dv^k_n&=&(D_i(a^{ij}_{kr}v^r_{nx^j}+\bar{b}^i_{kr}v^r_n+\bar{f}^{ik}_n)
+b^i_{kr}v^r_{nx^i}+c_{kr}v^r_n+f^k_n-a^{ij}_{kr}u^r_{x^j}\zeta_{nx^i})\,dt\nonumber\\
&&+(\sigma^{ik}_{kr,m}v^r_{nx^i}+\nu_{kr,m}v^r_n+g^k_m)\,dw^m_t,\label{eqn
5.2.1}
\end{eqnarray}
where
$$
f^k_n:=-(\bar{b}^i_{kr}u^r+\bar{f}^{ik}+b^i_{kr}u^r)\zeta_{nx^i}+f^k\zeta_n,
$$
$$
\bar{f}^{ik}_n:=-a^{ij}_{kr}u^r\zeta_{nx^j}+\bar{f}^{ik}\zeta_n,\qquad
g^k_n=-\sigma^{ik}u\zeta_{nx^i}+g^k\zeta_n.
$$
Also, we note that $\zeta_0u\in \cH^1_2(T)$ and
$\|\psi^{-1}\zeta_0u\|_{\bH^1_{2,\theta}(\cO,T)}\sim
\|\zeta_0u\|_{\bH^1_2(T)}$.  By Theorem \ref{thm 1} and Lemma
\ref{lemma a priori}, we have
\begin{equation}
                 \label{new 1}
\|\psi^{-1}u\|^2_{\bH^1_{2,\theta}(\cO,T)}\leq\sum_{n=0}^{N_0}
\|\psi^{-1}v_n\|^2_{\bH^1_{2,\theta}(\cO,T)}
\end{equation}
\begin{equation}
                          \label{eqn 8.13.1}
\leq N\sum_{n=0}^{N_0}( \|\bar{f}^i_n\|^2_{\bL_{2,\theta}(\cO,T)}+\|\psi
f_n\|^2_{\bH^{-1}_{2,\theta}(\cO,T)} + \|\psi
a^{ij}_{kr}u^r_{x}\zeta_{nx}\|^2_{\bH^{-1}_{2,\theta}(\cO,T)}
+\|g_n\|^2_{\bL_{2,\theta}(\cO,T,\ell_2)}).
\end{equation}
Actually relations like (\ref{new 1}) hold even if $N_0=\infty$ and this is why 
the theorem is true even when $\cO$ is not bounded.

Since  $a^{ij}$ is only measurable,
 at most we get
\begin{eqnarray*} \sum_n\|\psi
a^{ij}_{kr}u^r_{x^j}\zeta_{nx^i}\|^2_{\bH^{-1}_{2,\theta}(\cO,T)}\leq
\sum_n\|\psi
a^{ij}_{kr}u^r_{x^j}\zeta_{nx^i}\|^2_{\bL_{2,\theta}(\cO,T)} \leq
N\|u_x\|^2_{\bL_{2,\theta}(\cO,T)}\leq
N\|\psi^{-1}u\|^2_{\bH^1_{2,\theta}(\cO,T)}
\end{eqnarray*}
and consequently (\ref{eqn 8.13.1}) only  leads us to the useless
inequality
$$
\|\psi^{-1}u\|^2_{\bH^1_{2,\theta}(\cO,T)}\leq
N\|\psi^{-1}u\|^2_{\bH^1_{2,\theta}(\cO,T)}+ \cdot\cdot\cdot.
$$
Hence, to avoid  estimating the norm $\|\psi
a^{ij}_{kr}u^r_{x^j}\zeta_{nx^i}\|_{\bH^{-1}_{2,\theta}(\cO,T)}$ we
proceed as  in \cite{Kim}.  We note that for each $k$ we have
$$
\psi^{-1}a^{ij}_{kr}u^r_{x^j}\zeta_{nx^i}\in
\psi^{-1}\bL_{2,\theta}(\cO,T).
$$
Thus, by Theorem 2.9 in \cite{Kim03}, for each $k$ \;the solution
$\bar{v}^k_n \in
 \frH^2_{2,\theta}(\cO,T)$ of the single equation
$$
dv=(\Delta v -\psi^{-1}a^{ij}_{kr}u^r_{x^j}\zeta_{nx^i} )dt, \quad
v(0)=0
$$
satisfies
\begin{equation}
                                                                \label{7.21.02}
\|\bar{v}^k_n\|_{\frH^2_{2,\theta}(\cO,T)}\leq N
\|a^{ij}_{kr}u^r_{x^j}\zeta_{nx^i}\|_{\bL_{2,\theta}(\cO,T)}\leq
N\|u_{x}\zeta_{nx}\|_{\bL_{2,\theta}(\cO,T)}
\end{equation}
and, by Lemma \ref{lemma 15.05}, for each $t\leq T$ we have
\begin{equation}
                                     \label{eqn 08.01.1}
\|\bar{v}^k_n\|^2_{\bH^1_{2,\theta}(\cO,t)}\leq N
t\|\bar{v}^k_n\|^2_{\frH^2_{2,\theta}(\cO,t)}\leq Nt
\|u_{x}\zeta_{nx}\|^2_{\bL_{2,\theta}(\cO,t)},
\end{equation}
where $N$ is independent of $T$ since we assume $T\leq 1$. Now, we
denote $\bar{u}^k_n:=\bar{v}^k_n\psi\bar{\zeta}_n$ and
$\bar{u}_n=(\bar{u}^1,\cdot,\bar{u}^{d_1}_n)$. Then $\bar{u}_n$
satisfies
$$
d\bar{u}^k_n=(\Delta \bar{u}^k_n + \hat{f}^k_n
-a^{ij}_{kr}u^r_{x^j}\zeta_{nx^i}) \,dt, \quad \bar{u}^k_n(0)=0,
$$
where $\hat{f}^k_n= -2\bar{v}^k_{nx^i}(\bar{\zeta}_n
\psi)_{x^i}-\bar{v}^k_n \Delta(\bar{\zeta}_n \psi)$. Finally, as we
  denote $u_n:=v_n-\bar{u}_n$, we find that
 $u_n$ satisfies
$$                                                            \label{04.02.01}
du^k_{n}=( D_{i}(a^{ij}_{kr}u^r_{nx^j}+
\bar{b}^iu_n+\bar{F}^{ik}_n)+ b^i_{kr}u^r_{nx^i}+ c_{kr} u^r_n
+F^k_n)\,dt
$$
\begin{equation}
                                                            \label{04.05.04}
+
(\sigma^{i}_{kr,m}u^r_{nx^i}+\nu_{kr,m}u^r_n+G^{k}_{n,m})\,dw^{m}_{t},
\end{equation}
where
$$
\bar{F}^{ik}_n=\bar{f}^{ik}_n+(a^{ij}_{kr}-\delta^{ij}\delta^{kr})\bar{u}^r_{nx^j}+\bar{b}^i_{kr}\bar{u}^r_n,
$$
$$
F^k_n=f^k_n+\hat{f}^k_n+b^i_{kr}\bar{u}^r_{nx^i}+c_{kr}\bar{u}^r_n,
\quad
 G^{k}_n=\sigma^{i}_{kr}\bar{u}^r_{nx^i}+\nu_{kr}\bar{u}^r_n+g^k_n.
$$
Then, by Lemmas \ref{lemma a priori}, for any $n\geq 1$ and $t\leq
T$ we have
\begin{equation}
                         \label{eqn 8.14.1}
\|\psi^{-1}u_n\|^2_{\bH^1_{2,\theta}(\cO,t)}\leq N
\|\bar{F}^i_n\|^2_{\bL_{2,\theta}(\cO,t)} + N \|\psi
F_n\|^2_{\bH^{-1}_{2,\theta}(\cO,t)}
+N\|G_n\|^2_{\bL_{2,\theta}(\cO,t)}.
\end{equation}
Also, since  $u\zeta_0$ has compact support in $\cO$, (\ref{eqn
8.14.1}) holds for $n=0$ by Theorem \ref{thm 1}. As we recall that
$\psi b, \psi \bar{b}, \psi^2 c, \psi_x, \psi \psi_{xx},
(\bar{\zeta}_n\psi)_x, \psi \Delta(\bar{\zeta}_n \psi)$ are bounded,
 $\|\cdot\|_{H^{-1}_{2,\theta}(\cO)}\leq \|\cdot\|_{L_{2,\theta}(\cO)}$, and
$$
\psi^{-1}\bar{u}_n=\bar{\zeta}_n\bar{v}_n, \quad
\bar{u}_{nx}=\bar{\zeta}_n\psi \bar{v}_{nx}+\bar{v}_n
(\bar{\zeta}_n\psi)_x,
$$
we get
\begin{eqnarray*}
&&\|\psi
(\hat{f}^k_n+b^i_{kr}\bar{u}^r_{nx^i}+c_{kr}\bar{u}^r_n)\|_{\bH^{-1}_{2,\theta}(\cO,t)}\\
&\leq& N\left( \|\psi \bar{v}_{nx}\|_{\bL_{2,\theta}(\cO,t)}+\|
\bar{v}_{n}\|_{\bL_{2,\theta}(\cO,t)}+\|\bar{u}_{nx}\|_{\bL_{2,\theta}(\cO,t)}
+\|\psi^{-1}\bar{u}_n\|^2_{\bL_{2,\theta}(\cO,t)}\right)\\
&\leq& N \left(\|\psi \bar{v}_{nx}\|_{\bL_{2,\theta}(\cO,t)}+\|
\bar{v}_{n}\|_{\bL_{2,\theta}(\cO,t)}\right)\\
&\leq& N \|\bar{v}_n\|_{\bH^1_{2,\theta}(\cO,t)}
\end{eqnarray*}
and it leads to
$$
\|\psi F_n\|^2_{\bH^{-1}_{2,\theta}(\cO,t)}\leq N\|\psi
f_n\|^2_{\bH^{-1}_{2,\theta}(\cO,t)}+N\|\bar{v}_n\|^2_{\bH^1_{2,\theta}(\cO,t)}.
$$
Also, by (\ref{eqn 08.01.1}) we have
$$
\|\bar{v}_n\|^2_{\bH^{1}_{2,\theta}(\cO,t)}\leq
Nt\|u_{x}\zeta_{nx}\|_{\bL_{2,\theta}(\cO,t)}
$$
and consequently
\begin{eqnarray*}
&&\sum_n\|\psi F_n\|^2_{\bH^{-1}_{2,\theta}(\cO,t)}\\
&\leq& N \sum_n \left(\|\psi f_n\|^2_{\bH^{-1}_{2,\theta}(\cO,t)}+
t\|u_x\zeta_{nx}\|^2_{\bL_{p,\theta}(\cO,t)}\right)\\
 &\leq& N
\|u\|^2_{\bL_{2,\theta}(\cO,t)}+Nt\|u_x\|^2_{\bL_{2,\theta}(\cO,t)}
+N\|\bar{f}\|^2_{\bL_{2,\theta}(\cO,t)}+N\|\psi
f\|^2_{\bH^{-1}_{2,\theta}(\cO,t)}.
\end{eqnarray*}
The sums
$$
\sum_n \|\bar{F}^i_n\|^2_{\bL_{2,\theta}(\cO,t)}, \quad \sum_n
\|G_n\|^2_{\bL_{2,\theta}(\cO,t)}.
$$
are estimated similarly. Then for each $t\leq T$ one gets
\begin{eqnarray*}
\|\psi^{-1}u\|^2_{\bH^1_{2,\theta}(\cO,t)}&\leq& N \sum_n
\|\psi^{-1}v_n\|^2_{\bH^1_{p,\theta}(\cO,t)}\\
&\leq& N \|\bar{f}\|^2_{\bL_{2,\theta}(\cO,T)}+ N \|\psi
f\|^2_{\bH^{-1}_{2,\theta}(\cO,T)}+
N\|g\|^2_{\bL_{2,\theta}(\cO,T)}\\
&&+ N \|u\|^2_{\bL_{2,\theta}(\cO,t)}+ N\cdot
t\|\psi^{-1}u\|^2_{\bH^1_{2,\theta}(\cO,t)}.
\end{eqnarray*}
Now, we choose $\varepsilon \in (0,1]$ such that for $t\leq T\leq
\varepsilon$
$$
N\cdot t\|u_x\|^2_{\bL_{2,\theta}(\cO,t)}\leq 1/2
\|\psi^{-1}u\|^2_{\bH^1_{2,\theta}(\cO,t)}.
$$
Then, by Lemma \ref{lemma 15.05}, for each $t\leq T$ we obtain
\begin{eqnarray}
\|u\|^2_{\frH^1_{2,\theta}(\cO,t)}&\leq& N
\int^t_0\|u\|^2_{\frH^1_{2,\theta}(\cO,s)}\,ds + N
\|\bar{f}\|^2_{\bL_{2,\theta}(\cO,T)}\nonumber\\
&&+ N \|\psi f\|^2_{\bH^{-1}_{2,\theta}(\cO,T)}+
N\|g\|^2_{\bL_{2,\theta}(\cO,T)}.
\end{eqnarray}
This and Gronwall's inequality lead to the a priori estimate for
$T\leq \varepsilon$.
\end{proof}

 For the case $T\geq \varepsilon$ we need the following lemma,
 which is  proved in \cite{Kim} for $d_1=1$.

\begin{lemma}
                      \label{lemma 5.2}
Let $d-1<\theta<d+1+p$,\; $t_0\leq T$, and  $u\in
\frH^{\gamma+2}_{p,\theta}(\cO,t_0)$ satisfy
 $$
 du^k(t)=f^k(t)dt + g^k_m(t) dw^m_t, \quad u(0)=0.
 $$
 Then there exists a unique
$\tilde{u}\in \frH^{\gamma+2}_{p,\theta}(\cO,T)$ such that
$\tilde{u}(t)=u(t)$ for $t\leq t_0 (a.s)$ and on $(0,T)$
\begin{equation}
                                                                   \label{7.10.12}
d\tilde{u}^k=(\Delta \tilde{u}^k(t)+\tilde{f}^k(t))dt+g^k I_{t\leq
t_0} dw^m_t,
\end{equation}
where $\tilde{f}=(f^k(t)-\Delta u^k(t))I_{t\leq t_0}$. Furthermore,
we have
\begin{equation}
                                                                 \label{10.8.11}
\|\tilde{u}\|_{\frH^{\gamma+2}_{p,\theta}(\cO,T)}\leq
N\|u\|_{\frH^{\gamma+2}_{p,\theta}(\cO,t_0)},
\end{equation}
where $N$ is independent of $u$ and $t_0$.

\end{lemma}

\begin{proof}
We note that for each $k$, $\tilde{f}^k\in
\psi^{-1}\bH^{\gamma}_{p,\theta}(\cO,T)$ and $ g^kI_{t\leq t_0}\in
\bH^{\gamma+1}_{p,\theta}(\cO,T)$. Thus, by Theorem  2.9 in
\cite{Kim03},  equation (\ref{7.10.12}) has a unique (real-valued)
solution $\tilde{u}^k\in \frH^{\gamma+2}_{p,\theta}(\cO,T)$ and we
have
\begin{equation}
                                                                 \label{10.8.11}
\|\tilde{u}^k\|_{\frH^{\gamma+2}_{p,\theta}(\cO,T)}\leq
N\|u^k\|_{\frH^{\gamma+2}_{p,\theta}(\cO,t_0)}.
\end{equation}
We define
$\tilde{u}=(\tilde{u}^1,\tilde{u}^2,\cdots,\tilde{u}^{d_1})$. To
show $\tilde{u}(t)=u(t) $ for $t\leq t_0$ we notice that, for $t\leq
t_0$, the function
 $v^k(t)=\tilde{u}^k(t)-u^k(t)$ satisfies the equation
$$
dv^k(t)=\Delta v^k \, dt, \quad v(0,\cdot)=0.
$$
Thus, by Theorem 2.9 in \cite{Kim03},  $v^k(t) =0$ for $t\leq t_0 \,
(a.e)$. The lemma is proved.
\end{proof}
We finish the proof of Theorem \ref{main theorem on domain}.

{\bf Proof of Theorem \ref{main theorem on domain}} \;\; As usual, we only prove that 
estimate (\ref{main a priori}) holds given that a solution $u$ already exists. For simplicity, we assume $u_0=0$.
See the proof of Theorem 5.1 in \cite{Kr99} for the general case.

Take an
integer $M\geq 2$ such that $T/M\leq \varepsilon$ and we denote
$t_n=Tn/M$. Assume that, for $n=1,2,...,M-1$, we have the estimate
(\ref{main a priori})
 with $t_n$ in place of $T$ (and $N$ depending only on $d,d_1,\theta, \delta, K, L$ and $T$).
 We use the induction on $n$.

 Let $u_n \in \frH^1_{2,\theta}$ be the
continuation of $u$ on $[t_n, T]$, which exists by Lemma \ref{lemma
5.2} with
 $\gamma=-1$ and $t_0=t_n$. As we denote $v_n:=u-u_n$, we have $v_n(t)=0$ for $t\leq
 t_n
 (a.s)$ and, for any $t\in [t_n, T]$ and $\phi \in
 C^{\infty}_0(\cO)$,
\begin{eqnarray*}
(v^k_n(t), \phi)&=&-\int^t_{t_n}( a^{ij}_{kr}v^r_{n x^j}+
\bar{b}^i_{kr} v^r_n+ \bar{f}^{ik}_n, \phi_{x^i})(s) ds
+\int^t_{t_n} (b^i_{kr}v^r_{n x^i}+c_{kr}v^r_n+f_n, \phi)(s) ds\\
&& + \int^t_{t_n}(\sigma^{i}_{kr,m}u^r_{nx^i}+\nu_{kr,m}
v^r_n+g^k_{n,m}, \phi)(s)dw^m_s,
\end{eqnarray*}
 where
$$
\bar{f}^{ik}_n:=(a^{ij}_{kr}-\delta^{ij}\delta^{kr})u^r_{nx^j}+\bar{b}^i_{kr}
u^r_n+ \bar{f}^{ik}, \quad f^k_n=b^i_{kr} u^r_{n
x^i}+c_{kr}u^r_m+f^{k},
$$
$$
g^k_n:=\sigma^{i}_{kr}u^r_{nx^i}+\nu_{kr} u^r_n+g^k.
$$
Next, instead of random processes on $[0, T]$ we consider processes
given on $[t_n, T]$ and introduces spaces
$\frH^{\gamma}_{p,\theta}(\cO,[t_n,T])$,
 $\bL_{p,\theta}(\cO,[t_n,t])$, $\bH^{\gamma}_{p,\theta}(\cO,[t_n,T])$ in a natural
 way. Then we get a counterpart of the previous result and conclude
that
\begin{eqnarray*}
&&\bE\int^{t_{n+1}}_{t_n}\|\psi^{-1}(u-u_n)(s)\|^2_{H^1_{2,\theta}(\cO)}
ds\\
&\leq& N \bE
\int^{t_{n+1}}_{t_n}(\|\bar{f}^i_n(s)\|^2_{L_{2,\theta}(\cO)}+\|\psi
f_n(s)\|^2_{H^{-1}_{2,\theta}(\cO)}+
\|g_n(s)\|^2_{L_{2,\theta}(\cO)}) ds.
\end{eqnarray*}
Thus, by the induction hypothesis we get
\begin{eqnarray*}
&&\bE \int^{t_{n+1}}_{0}\|\psi^{-1}u(s)\|^2_{H^1_{2,\theta}(\cO)} ds\\
&\leq& N \bE \int^T_0 \|\psi^{-1}u_n(s)\|^2_{H^1_{2,\theta}(\cO)} ds +
N \bE
\int^{t_{n+1}}_{t_n}\|\psi^{-1}(u-u_n)(s)\|^2_{H^1_{2,\theta}(\cO)}
ds\\
&\leq& N (\|\bar{f}^i\|^2_{\bL_{2,\theta}(\cO, t_{n+1})}+ \|\psi
f\|^2_{\bH^{-1}_{2,\theta}(\cO,t_{n+1})}+\|g\|^2_{\bL_{2,\theta}(\cO,t_{n+1},\ell_2)}).
\end{eqnarray*}
We see that the induction goes through and thus the theorem is
proved.

\end{document}